\documentstyle[amssymb,12pt]{article}

\newlength{\minitwocolumn}
\newcommand{\Z}{{\Bbb Z}} 
\newcommand{\C}{{\Bbb C}} 

\newcommand{\la}{\lambda}

\newcommand{\nn}{{\nonumber}}
\newcommand{\eqref}[1]{(\ref{#1})}
\newcommand{\bea}{\begin{eqnarray}}
\newcommand{\ena}{\end{eqnarray}}
\newcommand{\beit}{\begin{itemize}}
\newcommand{\enit}{\end{itemize}}

\newcommand{\be}{\begin{eqnarray*}}
\newcommand{\en}{\end{eqnarray*}}
\newcommand{\lb}[1]{\label{#1}}
\newcommand{\End}{{\rm End}}

\newcommand{\id}{\hbox{id}}

\newcommand{\BW}[5]
{\left(\left.\matrix{#1 & #2 \cr #3 & #4 \cr}
\right| #5\right)}

\def\infq4p#1{{(#1;q^4,p)_\infty}}

\newcommand{\vep}{\varepsilon}

\font\teneufm=eufm10
\font\seveneufm=eufm7
\font\fiveeufm=eufm5
\newfam\eufmfam
\textfont\eufmfam=\teneufm
\scriptfont\eufmfam=\seveneufm
\scriptscriptfont\eufmfam=\fiveeufm

\let\goth\frak
\newcommand{\slth}{\widehat{\goth{sl}}_2}

\makeatletter
\@addtoreset{equation}{section}
\makeatother

\newtheorem{thm}{Theorem}[section]
\newtheorem{prop}[thm]{Proposition}

\newtheorem{cor}[thm]{Corollary}

\newcommand{\vt}[3]{\vartheta_{#1}\left(\left. \frac{#2}{r} \right| {#3} \right)}
\newcommand{\vtf}[3]{\vartheta_{#1}\left(\left. {#2} \right| {#3} \right)}
\newcommand{\vth}[2]{\vartheta_{#1}\left(\left. \frac{#2}{2r} \right| \frac{\tau}{2} \right)}

\newcommand{\sn}{{\rm sn}}
\newcommand{\cn}{{\rm cn}}
\newcommand{\dn}{{\rm dn}}

\newcommand{\tp}{{\tilde{p}}}

\newcommand{\baj}{{\bar{j}}}
\newcommand{\bao}{{\bar{1}}}
\newcommand{\bat}{{\bar{2}}}

\begin{document}

\begin{center}
{\Large\bf Fusion of Baxter's Elliptic $R$-matrix\\
and the Vertex-Face Correspondence\footnote{To appear in the proceedings of the workshop ``Solvable Lattice Models 2004", July 20--23, 2004, {\it RIMS Koukyuroku}, RIMS, Kyoto University.}}

\vspace{15mm}
{\large  Hitoshi Konno }

\vspace{3mm}
{\it 
Department of Mathematics, Faculty of Integrated Arts $\&$ Sciences, \\
Hiroshima University, Higashi Hiroshima 739-8521, Japan\\
E-mail:konno@mis.hiroshima-u.ac.jp}

\vspace{20mm}
{\bf Abstract}
\end{center}

\noindent
The matrix elements of the $2\times 2$ fusion of Baxter's elliptic $R$-matrix, $R^{(2,2)}(u)$, are given explicitly. Based on a note by Jimbo, we give a formula which show that $R^{(2,2)}(u)$ is gauge equivalent to Fateev's $R$-matrix for the 21-vertex model. 
Then the crossing symmetry formula for $R^{(2,2)}(u)$ is derived. 
We also consider the fusion of the vertex-face correspondence relation and derive 
a crossing symmetry relation between the fusion of 
the intertwining vectors and their dual vectors.

\section{Notations}
Let ${p}=e^{-\frac{\pi K'}{K}}$, $q=-e^{-\frac{\pi \lambda}{2K}}$ and $\zeta=e^{-\frac{\pi \la u}{2K}}$.
We introduce $x$, $\tau$ and $r$ by $x=-q$, 
$\tau=\frac{2iK}{K'}$ and $r=\frac{K'}{\la}$. Then $p=e^{-\frac{2\pi i}{ \tau}}=x^{2r}$. 
Through this paper, we assume ${\rm Im}\tau >0$.
Let $\tp=e^{2\pi i \tau}=e^{-\pi \frac{I'}{I}}$, where $I=\frac{K'}{2},\ 
I'=2K$. We use the theta functions 
\be
&&\vartheta_1(u|\tau)=2\tilde{p}^{1/8}(\tilde{p};\tilde{p})_\infty\sin\pi u
\prod_{n=1}^\infty(1-2\tilde{p}^n\cos2\pi u+\tilde{p}^{2n}),\\
&&\vartheta_0(u|\tau)=-ie^{\pi i(u+\tau/4)}\vartheta_1\left(u+\frac{\tau}{2}\Big|\tau\right),\\
&&\vartheta_2(u|\tau)=\vartheta_1\left(u+\frac{1}{2}\Big|\tau\right),\\
&&\vartheta_3(u|\tau)=e^{\pi i(u+\tau/4)}\vartheta_1\left(u+\frac{\tau+1}{2}\Big|\tau\right)
\en
and Jacobi's elliptic functions
\be
&&\sn\la u=\frac{\vtf{3}{0}{\tau}\vtf{1}{\frac{\la u}{2I}}{\tau}}{\vtf{2}{0}{\tau}\vtf{0}{\frac{\la u}{2I}}{\tau}}
=\frac{\vtf{3}{0}{\tau}\vtf{1}{\frac{u}{r}}{\tau}}{\vtf{2}{0}{\tau}\vtf{0}{\frac{u}{r}}{\tau}},\\
&&\cn\la u=\frac{\vtf{0}{0}{\tau}\vtf{2}{\frac{\la u}{2I}}{\tau}}{\vtf{2}{0}{\tau}\vtf{0}{\frac{\la u}{2I}}{\tau}}
=\frac{\vtf{0}{0}{\tau}\vtf{2}{\frac{u}{r}}{\tau}}{\vtf{2}{0}{\tau}\vtf{0}{\frac{u}{r}}{\tau}},\\
&&\dn\la u=\frac{\vtf{0}{0}{\tau}\vtf{3}{\frac{\la u}{2I}}{\tau}}{\vtf{3}{0}{\tau}\vtf{0}{\frac{\la u}{2I}}{\tau}}
=\frac{\vtf{0}{0}{\tau}\vtf{3}{\frac{u}{r}}{\tau}}{\vtf{3}{0}{\tau}\vtf{0}{\frac{u}{r}}{\tau}}.
\en
We also use the symbol $[u]$ defined by 
\be
&&[u]=x^{\frac{u^2}{r}-u}\Theta_{x^{2r}}(x^{2u})=C\vt{1}{u}{\tau}, \quad C=x^{-\frac{r}{4}}e^{-\frac{\pi i}{4}}\tau^{\frac{1}{2}}
\en
and abbreviation 
\be
&&\vth{1,2}{u}=\vth{1}{u}\vth{2}{u}=\vtf{0}{0}{\tau}\vt{1}{u}{\tau},
\en
etc..

\section{Fusion of Baxter's $R$-matrix}
Baxter's elliptic $R$-matrix is given by\cite{Baxter}
\bea
&&{R}(u)=
{R}_0(u)\left(\matrix{a(u)&&&d(u)\cr
&b(u)&c(u)&\cr
&c(u)&b(u)&\cr 
d(u)&&&a(u)\cr}\right),
\ena
where
\bea
R_0(u)&=&z^{-\frac{r-1}{2r}}\frac{(px^2z;x^4,p)_\infty(x^2z;x^4,p)_\infty
(p/z;x^4,p)_\infty(x^4/z;x^4,p)_\infty}{(px^2/z;x^4,p)_\infty
(x^2/z;x^4,p)_\infty(pz;x^4,p)_\infty(x^4z;x^4,p)_\infty},\lb{evR}\\
a(u)&=&\frac{\vtf{2}{\frac{1}{2r}}{\frac{\tau}{2}}\vtf{2}{\frac{  u}{2r}}{\frac{\tau}{2}}}{\vtf{2}{0}{\frac{\tau}{2}}\vtf{2}{\frac{1+u}{2r}}{\frac{\tau}{2}}},\qquad
b(u)=\frac{\vtf{2}{\frac{1 }{2r}}{\frac{\tau}{2}}\vtf{1}{\frac{  u}{2r}}{\frac{\tau}{2}}}
{\vtf{2}{0}{\frac{\tau}{2}}\vtf{1}{\frac{1+u}{2r}}{\frac{\tau}{2}}},\\
c(u)&=&\frac{\vtf{1}{\frac{ 1}{2r}}{\frac{\tau}{2}}\vtf{2}{\frac{  u}{2r}}{\frac{\tau}{2}}}
{\vtf{2}{0}{\frac{\tau}{2}}\vtf{1}{\frac{1+u}{2r}}{\frac{\tau}{2}}},\qquad
d(u)=-\frac{\vtf{1}{\frac{ 1}{2r}}{\frac{\tau}{2}}\vtf{1}{\frac{  u}{2r}}{\frac{\tau}{2}}}
{\vtf{2}{0}{\frac{\tau}{2}}\vtf{2}{\frac{1+u}{2r}}{\frac{\tau}{2}}}
\ena
with $z=\zeta^2=x^{2u}$.
Let $V=\C v_{\vep_1}\oplus \C v_{\vep_2},\ \vep_1,\vep_2=+,-$. We regard $R(u)\in \End(V\otimes V)$.  
The $R$-matrix \eqref{evR} satisfies
\bea
&&R(u)PR(u)P=\id,\\
&&R(-u-1)=(\sigma^{y }\otimes 1)^{-1}\ (PR(u)P)^{t_1} \  \sigma^y\otimes 1,\\
&&R(0)=P, \qquad \lim_{u\to -1}R(u)=P-\id.
\ena
Here ${^{t_1}}$ denotes the transposition with respect to the first vector space in the tensor product and $P(\vep_1\otimes \vep_2)=\vep_2\otimes \vep_1$.

Fusion of $R(u)$ was considered systematically in \cite{DJKMO}. 
Let $V^{(2)}$ be the space of the symmetric tensors in $V\otimes V$      
spanned by $v^{(2)}_{2}\equiv v_+\otimes v_+,\ v^{(2)}_{0}\equiv 
\frac{1}{2}(v_+\otimes v_-+v_-\otimes v_+),\ v^{(2)}_{-2}\equiv v_-\otimes v_-$.The projection operator of the space $V\otimes V$ on  $V^{(2)}$ 
is given by $\Pi=\frac{1}{2}(P+\id)$. 
Let $V_1, V_2, V_{\bao}, V_{\bat}$ be the copies of V. 
Define 
\bea
&&R^{(2,1)}_{12, \bar{j}}(u)=\Pi_{12}R_{1\baj}(u+1)R_{2\baj}(u) \ \in {\rm End}(V^{(2)} \otimes V_{\baj}).
\lb{hfusion}
\ena
It follows that  
\bea
&&R^{(2,1)}_{12,\baj}(u)\Pi_{12}=R^{(2,1)}_{12,\baj}(u).
\ena
The 2$\times$2 fusion of the $R$-matrix is then given by
\bea
&&R^{(2,2)}(u)=\Pi_{\bao\bat}R^{(2,1)}_{12,\bat}(u)R^{(2,1)}_{12,\bao}(u-1) \ \in {\rm End}(V^{(2)} \otimes V^{(2)}).\lb{vhfusion}
\ena
This satisfies the Yang-Baxter equation on $V^{(2)} \otimes V^{(2)} \otimes V^{(2)}$. 

We calculate the matrix elements of $R^{(2,1)}(u)$ and $R^{(2,2)}(u)$ defined by  
\bea 
&&R^{(2,1)}(u) v^{(2)}_{\mu}\otimes v^{}_{\vep}=\sum_{\mu'=2,0,-2 \atop \vep'=+,-} R^{(2,1)}(u)^{\mu \vep}_{\mu' \vep'}\ v^{(2)}_{\mu'}\otimes v^{}_{\vep'},\\
&&R^{(2,2)}(u) v^{(2)}_{\mu_1}\otimes v^{(2)}_{\mu_2} =\sum_{\mu'_1,\mu'_2=2,0,-2} R^{(2,2)}(u)^{\mu_1 \mu_2}_{\mu'_1 \mu'_2}\ v^{(2)}_{\mu'_1}\otimes v^{(2)}_{\mu'_2}.
\ena
From \eqref{hfusion}, we have 
\bea
&&R^{(2,1)}(u)^{\mu\ \bar{\vep}}_{\mu'\ \bar{\vep}'}=\sum_{\vep_2',\bar{\vep}''=\pm 1}R(u+1)^{\mu-\vep_2\ \bar{\vep}''}_{\mu'-\vep_2'\ \bar{\vep}'}R(u)^{\vep_2\ \bar{\vep}}_{\vep_2'\ \bar{\vep}''},
\ena
where we set $\mu=\vep_1+\vep_2,\  \mu'=\vep_1'+\vep_2'$. Evaluating the summation explicitly, we obtain  
\begin{prop}
\be
&&R^{(2,1)}(u)=R^{(2,1)}_0(u)\left(\matrix{R^{+2+}_{+2+}&0&0&R^{\ 0-}_{+2+}&R^{-2+}_{+2+}&0\cr 
                 0&R^{+2-}_{+2-}&R^{\ 0+}_{+2-}&0&0&R^{-2-}_{+2-}\cr
                 0&R^{+2-}_{\ 0+}&R^{0+}_{0+}&0&0&R^{-2-}_{\ 0+}\cr
                R^{+2+}_{\ 0-}&0&0&R^{\ 0-}_{\ 0-}&R^{-2+}_{\ 0-}&0\cr
                R^{+2+}_{-2+}&0&0&R^{\ 0-}_{-2+}&R^{-2+}_{-2+}&0\cr
                0&R^{+2-}_{-2-}&R^{\ 0+}_{-2-}&0&0&R^{-2-}_{-2-}\cr}\right),
\en
where 
\be
&&R^{(2,1)}_0(u)=R_0(u+1)R_0(u)=-\frac{[u+1]}{[u]},\\
&&R(u)^{+2+}_{+2+}=R(u)^{-2-}_{-2-}=\frac{\vtf{2}{\frac{1}{2r}}{\frac{\tau}{2}}^2\vtf{2}{\frac{u}{2r}}{\frac{\tau}{2}}}{\vtf{2}{0}{\frac{\tau}{2}}^2\vtf{2}{\frac{u+2}{2r}}{\frac{\tau}{2}}},\\
&&R(u)^{\ 0-}_{+2+}=R(u)^{\ 0+}_{-2-}
=-\frac{\vtf{1}{\frac{1}{2r}}{\frac{\tau}{2}}\vtf{2}{\frac{1}{2r}}{\frac{\tau}{2}}
\vtf{1}{\frac{u}{2r}}{\frac{\tau}{2}}}{\vtf{2}{0}{\frac{\tau}{2}}^2\vtf{2}{\frac{u+2}{2r}}{\frac{\tau}{2}}},\\
&&R(u)^{-2+}_{+2+}=R(u)^{+2-}_{-2-}=-\frac{\vtf{1}{\frac{1}{2r}}{\frac{\tau}{2}}^2\vtf{2}{\frac{u}{2r}}{\frac{\tau}{2}}}{\vtf{2}{0}{\frac{\tau}{2}}^2\vtf{2}{\frac{u+2}{2r}}{\frac{\tau}{2}}},\\
&&R(u)^{+2-}_{+2-}=R(u)^{-2+}_{-2+}=\frac{\vtf{2}{\frac{1}{2r}}{\frac{\tau}{2}}^2\vtf{1}{\frac{u}{2r}}{\frac{\tau}{2}}}{\vtf{2}{0}{\frac{\tau}{2}}^2\vtf{1}{\frac{u+2}{2r}}{\frac{\tau}{2}}},\\
&&R(u)^{\ 0+}_{+2-}=R(u)^{\ 0-}_{-2+}=\frac{\vtf{1}{\frac{1}{2r}}{\frac{\tau}{2}}\vtf{2}{\frac{1}{2r}}{\frac{\tau}{2}}
\vtf{2}{\frac{u}{2r}}{\frac{\tau}{2}}}{\vtf{2}{0}{\frac{\tau}{2}}^2\vtf{1}{\frac{u+2}{2r}}{\frac{\tau}{2}}},\\
&&R(u)^{-2-}_{+2-}=R(u)^{+2+}_{-2+}=-
\frac{\vtf{1}{\frac{1}{2r}}{\frac{\tau}{2}}^2\vtf{1}{\frac{u}{2r}}{\frac{\tau}{2}}}{\vtf{2}{0}{\frac{\tau}{2}}^2\vtf{1}{\frac{u+2}{2r}}{\frac{\tau}{2}}},\\
&&R(u)^{+2-}_{0\ +}=R(u)^{-2+}_{0\ -}=
\frac{\vtf{1}{\frac{1}{r}}{\frac{\tau}{2}}\vtf{2}{\frac{u+1}{2r}}{\frac{\tau}{2}}^2}
{\vtf{2}{0}{\frac{\tau}{2}}
\vtf{1}{\frac{u+2}{2r}}{\frac{\tau}{2}}\vtf{2}{\frac{u+2}{2r}}{\frac{\tau}{2}}},\\
&&R(u)^{+2+}_{0\ -}=R(u)^{-2-}_{0\ +}=-
\frac{\vtf{1}{\frac{1}{r}}{\frac{\tau}{2}}\vtf{1}{\frac{u+1}{2r}}{\frac{\tau}{2}}^2}
{\vtf{2}{0}{\frac{\tau}{2}}
\vtf{1}{\frac{u+2}{2r}}{\frac{\tau}{2}}\vtf{2}{\frac{u+2}{2r}}{\frac{\tau}{2}}},\\
&&R(u)^{0\ +}_{0\ +}=R(u)^{0\ -}_{0\ -}=
\frac{\vtf{2}{\frac{1}{r}}{\frac{\tau}{2}}\vth{1,2}{u+1}}
{\vtf{2}{0}{\frac{\tau}{2}}\vth{1,2}{u+2}}.
\en
\end{prop}

Similarly, from \eqref{vhfusion}, we obtain  
\be
&&R^{(2,2)}(u)^{\mu_1\ \mu_2}_{\mu_1'\ \mu_2'}=\sum_{\mu''=0, \pm2 \atop \bar{\vep}_1'=\pm1}R^{(2,1)}(u)^{ \mu''\ {\mu}_2-\bar{\vep}_1}_{\mu_1'\ {\mu}_2'-\bar{\vep}_1'}R^{(2,1)}(u-1)^{ \mu_1\ \bar{\vep}_1}_{\mu''\ \bar{\vep}_1'},
\en
where $\mu_1=\vep_1+\vep_2,\ \mu_1'=\vep_1'+\vep_2'$ and 
$\mu_2=\bar{\vep}_{{1}}+\bar{\vep}_{{2}},\ \mu_2'=\bar{\vep}_{{1}}'
+\bar{\vep}_{{2}}'$. 
\begin{prop}
\be
&&R^{(2,2)}(u)=R^{(2,2)}_0(u)\left(\matrix{
G&0&A&0&B&0&A&0&H\cr 
0&F&0&C&0&C^*&0&D&0\cr
A^*&0&G^*&0&B^*&0&H^*&0&A^*\cr 
0&C&0&F&0&D&0&C^*&0\cr
I&0&I^*&0&E&0&I^*&0&I\cr 
0&C^*&0&D&0&F&0&C&0\cr
A^*&0&H^*&0&B^*&0&G^*&0&A^*\cr 
0&D&0&C^*&0&C&0&F&0\cr
H&0&A&0&B&0&A&0&G\cr 
}\right),
\en
where
\be
R^{(2,2)}_0(u)&=&R^{(2,1)}_0(u-1)R^{(2,1)}_0(u)=\frac{[u+1]}{[u-1]},\\
A&=&R(u)^{+2-2}_{+2+2}=R(u)^{-2+2}_{-2-2}=R(u)^{-2+2}_{+2+2}=R(u)^{+2-2}_{-2-2}\nn\\
&=&-
\frac{
\vtf{1}{\frac{1}{r}}{\frac{\tau}{2}}
\vtf{2}{\frac{u}{2r}}{\frac{\tau}{2}}
\vtf{1,2}{\frac{1}{2r}}{\frac{\tau}{2}}
\vtf{1,2}{\frac{u}{2r}}{\frac{\tau}{2}}
}{
\vtf{2}{0}{\frac{\tau}{2}}^3
\vtf{2}{\frac{u+2}{2r}}{\frac{\tau}{2}}
\vtf{1,2}{\frac{u+1}{2r}}{\frac{\tau}{2}}
},\\
B&=&R(u)^{\ 0\ 0}_{+2+2}=R(u)^{\ 0\ 0}_{-2-2}=-
\frac{
\vtf{2}{\frac{1}{r}}{\frac{\tau}{2}}
\vtf{1}{\frac{u}{2r}}{\frac{\tau}{2}}
\vtf{1,2}{\frac{1}{2r}}{\frac{\tau}{2}}
\vtf{1,2}{\frac{u}{2r}}{\frac{\tau}{2}}
}{
\vtf{2}{0}{\frac{\tau}{2}}^3
\vtf{2}{\frac{u+2}{2r}}{\frac{\tau}{2}}
\vtf{1,2}{\frac{u+1}{2r}}{\frac{\tau}{2}}
},\\
C&=&R(u)^{+2\ 0}_{\ 0+2}=R(u)^{-2\ 0}_{\ 0-2}=R(u)^{\ 0+2}_{+2\ 0}=R(u)^{\ 0-2}_{-2\ 0}=
\frac{
\vtf{2}{\frac{u}{2r}}{\frac{\tau}{2}}^2 \vtf{1,2}{\frac{1}{r}}{\frac{\tau}{2}}
}
{
\vtf{2}{0}{\frac{\tau}{2}}^2 \vth{1,2}{u+2}
},\\
D&=&R(u)^{+2\ 0}_{-2\ 0 }=R(u)^{-2\ 0}_{+2\ 0 }=R(u)^{\ 0+2}_{\ 0-2}=R(u)^{\ 0 -2}_{\ 0+2}
=-\frac{
\vtf{1}{\frac{1}{r}}{\frac{\tau}{2}}^2 \vth{1,2}{u}
}
{
\vtf{2}{0}{\frac{\tau}{2}}^2 \vth{1,2}{u+2}
},\\
E&=&R(u)^{\ 0\ 0 }_{\ 0\ 0 }\nn\\
&=&\frac{
\vtf{1}{\frac{1}{r}}{\frac{\tau}{2}}
\vtf{1,2}{\frac{1}{2r}}{\frac{\tau}{2}}
\left(
\vtf{2}{\frac{u+1}{2r}}{\frac{\tau}{2}}^3
\vtf{2}{\frac{u-1}{2r}}{\frac{\tau}{2}}+
\vtf{1}{\frac{u+1}{2r}}{\frac{\tau}{2}}^3
\vtf{1}{\frac{u-1}{2r}}{\frac{\tau}{2}}
\right)
}{
\vtf{2}{0}{\frac{\tau}{2}}^3
\vtf{1,2}{\frac{u+2}{2r}}{\frac{\tau}{2}}
\vtf{1,2}{\frac{u+1}{2r}}{\frac{\tau}{2}}}\nn\\
&&+\frac{
\vtf{2}{\frac{1}{r}}{\frac{\tau}{2}}^2\vth{1,2}{u}
}{
\vtf{2}{0}{\frac{\tau}{2}}^2\vth{1,2}{u+2}
}
,\\
F&=&R(u)^{\ 0+2}_{\ 0 +2}=R(u)^{\ 0-2}_{\ 0 -2}=R(u)^{+2\ 0}_{+2\ 0}=R(u)^{-2\ 0}_{-2\ 0}
=\frac{
\vtf{2}{\frac{1}{r}}{\frac{\tau}{2}}^2 \vth{1,2}{u}
}
{
\vtf{2}{0}{\frac{\tau}{2}}^2 \vth{1,2}{u+2}
},
\en
\be
G&=&R(u)^{+2+2}_{+2+2}=R(u)^{-2-2}_{-2-2}\nn\\
&=&
\frac{
\vtf{2}{\frac{u}{2r}}{\frac{\tau}{2}}
\left(\vtf{1}{\frac{1}{2r}}{\frac{\tau}{2}}^4\vtf{1}{\frac{u-1}{2r}}{\frac{\tau}{2}}
\vtf{2}{\frac{u+1}{2r}}{\frac{\tau}{2}}+
\vtf{2}{\frac{1}{2r}}{\frac{\tau}{2}}^4\vtf{2}{\frac{u-1}{2r}}{\frac{\tau}{2}}
\vtf{1}{\frac{u+1}{2r}}{\frac{\tau}{2}}
\right)
}
{
\vtf{2}{0}{\frac{\tau}{2}}^4 \vtf{2}{\frac{u+2}{2r}}{\frac{\tau}{2}}
\vtf{1,2}{\frac{u+1}{2r}}{\frac{\tau}{2}}
}
,
\en
\be
H&=&R(u)^{+2+2}_{-2-2}=R(u)^{-2-2}_{+2+2}=
\frac{
\vtf{1}{\frac{1}{r}}{\frac{\tau}{2}}
\vtf{1}{\frac{u}{2r}}{\frac{\tau}{2}}^3
\vth{1,2}{1}
}{
\vtf{2}{0}{\frac{\tau}{2}}^3
\vtf{2}{\frac{u+2}{2r}}{\frac{\tau}{2}}
\vth{1,2}{u+1}
},
\en
\be
I&=&R(u)^{+2+2}_{\ 0\ 0 }=R(u)^{-2-2}_{\ 0\ 0}\nn\\
&=&-\frac{
\vtf{1}{\frac{1}{r}}{\frac{\tau}{2}}
\left(
\vtf{2}{\frac{1}{2r}}{\frac{\tau}{2}}^2
\vtf{1}{\frac{u+1}{2r}}{\frac{\tau}{2}}^3
\vtf{2}{\frac{u-1}{2r}}{\frac{\tau}{2}}+
\vtf{1}{\frac{1}{2r}}{\frac{\tau}{2}}^2
\vtf{2}{\frac{u+1}{2r}}{\frac{\tau}{2}}^3
\vtf{1}{\frac{u-1}{2r}}{\frac{\tau}{2}}
\right)
}{
\vtf{2}{0}{\frac{\tau}{2}}^3
\vtf{1,2}{\frac{u+2}{2r}}{\frac{\tau}{2}}
\vtf{1,2}{\frac{u+1}{2r}}{\frac{\tau}{2}}}\nn\\
&&-\frac{
\vtf{1}{\frac{u}{2r}}{\frac{\tau}{2}}^2\vtf{1,2}{\frac{1}{r}}{\frac{\tau}{2}}
}{
\vtf{2}{0}{\frac{\tau}{2}}^2\vth{1,2}{u+2}
}.
\en
The $*$-ed matrix element is obtained from a corresponding non-$*$-ed element by replacing 
the theta functions only depending on $u$ in the following rule.
\be
&&\vartheta_1 \to -\vartheta_2,\quad \vartheta_2 \to \vartheta_1. 
\en
\end{prop}
From this expression, one can easily see the following symmetries.
\be
&&R^{(2,2)}(u)^{\vep_1 \vep_2}_{\vep_1' \vep_2'}=R^{(2,2)}(u)^{\vep_2 \vep_1}_{\vep_2' \vep_1'}, \qquad\qquad ( P{\rm-invariance})\\
&&R^{(2,2)}(u)^{\vep_1 \vep_2}_{\vep_1' \vep_2'}=R^{(2,2)}(u)^{-\vep_1 -\vep_2}_{-\vep_1' -\vep_2'} \qquad\qquad ({\rm \Z_2-symmetry}).
\en 

In 1980, Fateev proposed the 21-vertex model as the spin one extension of Baxter's eight vertex model\cite{Fateev}. 
Solving the Yang-Baxter equation, he obtained the following $R$-matrix.
\be
&&R_{F}(u)=\tilde{F}(u)\left(\matrix{
s_1&0&0&0&\mu&0&0&0&\nu\cr
0&t&0&r&0&0&0&0&0\cr
0&0&T&0&0&0&R&0&0\cr
0&r&0&t&0&0&0&0&0\cr
\mu&0&0&0&s_2&0&0&0&\rho\cr
0&0&0&0&0&a&0&q&0\cr
0&0&R&0&0&0&T&0&0\cr
0&0&0&0&0&q&0&a&0\cr
\nu&0&0&0&\rho&0&0&0&s_3\cr    
}
\right)
\en
where
\be
&&s_1=\cn2\la+\frac{\sn\la\ \sn2\la}{\sn\la u\ \sn\la (u+1)},\\
&&s_2=\cn2\la+\dn2\la-1+\frac{\sn\la\ \sn2\la}{\sn\la u\ \sn\la( u+1)},\\
&&s_3=\dn2\la+\frac{\sn\la\ \sn2\la}{\sn\la u\ \sn\la( u+1)},\\
&&T=1,\qquad\qquad t=\cn2\la,\qquad\qquad a=\dn2\la,\\
&&r=\frac{\cn\la u\ \sn2\la}{\sn\la u},\qquad\qquad \mu=-\frac{\cn\la( u+1)\ \sn2\la}{\sn\la( u+1)},\\
&&R=\frac{ \sn2\la}{\sn\la u},\qquad\qquad \nu=-\frac{ \sn2\la}{\sn\la( u+1)},\\
&&q=\frac{\dn\la u\ \sn2\la}{\sn\la u},\qquad\qquad \rho=-\frac{\dn\la( u+1)\ \sn2\la}{\sn\la( u+1)}
\en
and $\tilde{F}(u)$ satisfies
\be
&&\tilde{F}(u)=\tilde{F}(-u-1),\\
&&\tilde{F}(u)\tilde{F}(-u)=\frac{\sn^2\la u}{\sn^2\la u-\sn^22\la}.
\en
$R_{F}(u)$ has the following symmetries.
\be
&&R_{F}(u)^{ij}_{kl}=R_{F}(u)^{ji}_{lk} \qquad\qquad ( P{\rm-invariance})\\
&&R_{F}(u)^{ij}_{kl}=R_{F}(u)_{ij}^{kl},\qquad\qquad (T{\rm-invariance})\\
&&R_{F}(u)^{ij}_{kl}=R_{F}(-u-1)_{il}^{kj}\qquad\qquad ({\rm Crossing\ symmetry}).
\en

We find 
\begin{prop}
\be
&&\tilde{F}(u)=R^{(2,2)}_0(u)\frac{\vt{0}{2}{\tau}\vt{1}{u}{\tau}}{\vtf{0}{0}{\tau}\vt{1}{u+2}{\tau}}.
\en
\end{prop}

Then the  following theorem is essentially due to Jimbo\cite{Jimbo}.
\begin{thm}
$R^{(2,2)}(u)$ is gauge equivalent to $R_{F}(u)$. Namely, 
\be
&&R_{F}(u)=U\otimes U \ R^{(2,2)}(u)\ (U\otimes U)^{-1},
\en
where
\be
&&U=\left(\matrix{1&0&0\cr 0&x&0\cr 0&0&y\cr}\right)\left(\matrix{\frac{1}{\sqrt{2}}&0&\frac{1}{\sqrt{2}}\cr 0&1&0\cr 
\frac{1}{\sqrt{2}}&0&-\frac{1}{\sqrt{2}}\cr}\right),
\en
and
\be
&&x^2={-\frac{1}{2}\frac{\vtf{0}{0}{\tau}\vt{3}{1}{\tau}}{\vt{0}{1}{\tau}\vtf{3}{0}{\tau}}},\quad 
y^2={-\frac{\vtf{2}{0}{\tau}\vt{3}{1}{\tau}}{\vt{2}{1}{\tau}\vtf{3}{0}{\tau}}}.
\en
\end{thm}
Combining the crossing symmetry of $R_{F}(u)$ and the $P$-invariance of $R^{(2,2)}(u)$, we find  
the following crossing symmetry formula for $R^{(2,2)}(u)$. 
\begin{cor}
\bea
&&R^{(2,2)}(-u-1)=Q^{-1}\otimes 1\ (P^{(2)}R^{(2,2)}(u)P^{(2)})^{t_1}\ Q\otimes 1,\lb{crossing}
\ena
where 
\be
&&Q=U^t U=\frac{1}{2}\left(\matrix{1+y^2&0&1-y^2\cr 0&x^2&0\cr 1-y^2&0&1+y^2\cr}\right)
\en
and $P^{(2)}$ is the permutation operator $P^{(2)}( v^{(2)}_{\vep_1}\otimes v^{(2)}_{\vep_2})= v^{(2)}_{\vep_2}\otimes v^{(2)}_{\vep_1}$.
\end{cor}

\noindent
{\it Remark :}\

\noindent
The crossing symmetry of the elliptic $R$-matrix is related to the dual module of the finite 
dimensional module of the elliptic algebra ${\cal A}_{q,p}(\slth)$, or the module of the Sklyanin algebra. See \cite{JM} for the case $U_q(\slth)$.  
Let $V_\zeta$ be the 3-dimensional module of ${\cal A}_{q,p}(\slth)$, and $V_\zeta^*$ its dual module. 
The above $Q$-matrix gives an isomorphism between $V_\zeta$ and $V_\zeta^*$. 

\section{The Vertex-Face Correspondence }
The vertex-face correspondence is a relationship between Baxter's $R$-matrix and the SOS face weight $W\BW{a_1}{a_2}{a_4}{a_3}{u}$ given by  
  \begin{eqnarray}
&&{W}\left(\left.
\begin{array}{cc}
n&n\pm 1\\
n\pm1&n\pm2
\end{array}\right|u\right)={R}_0(u),\nn\\
&&{W}\left(\left.
\begin{array}{cc}
n&n\pm 1\\
n\pm1&n
\end{array}\right|u\right)={R}_0(u)\frac{[n\mp u][1]}{[ n][ 1+u]},\lb{face}\\
&&{W}\left(\left.
\begin{array}{cc}
n&n\pm 1\\
n\mp 1&n
\end{array}\right|u\right)=
{R}_0(u)
\frac{[ n\pm 1][u] }{[ n][1+u]}.\nn
\end{eqnarray}
Let us consider the following vectors in $V$   
\bea
&&\psi(u)^a_b=\psi_+(u)^a_b\ v_+ + \psi_-(u)^a_b\ v_-,\qquad \\
&&{\psi}_+(u)^a_b=\vtf{0}{\frac{(a-b)u+a}{2r}}{\frac{\tau}{2}},\qquad \psi_-(u)=\vtf{3}{\frac{(a-b)u+a}{2r}}{\frac{\tau}{2}} \lb{intertwinvec}\nn
\ena 
with $|a-b|=1$. Baxter showed the following identity\cite{Baxter}.
\bea
&&\sum_{\vep_1',\vep_2'}
R(u-v)_{\vep_1 \vep_2}^{\vep_1' \vep_2'}\ 
\psi_{\vep_1'}(u)_{b}^{a}
\psi_{\vep_2'}(v)_{c}^{b}
=\sum_{b' \in {\mathbb{Z}}}
\psi_{\vep_2}(v)_{b'}^{a}
\psi_{\vep_1}(u)_{c}^{b'}
W\left(\left.
\begin{array}{cc}
a&b\\
b'&c
\end{array}\right|u-v\right).\lb{vertexface}
\ena

\subsection{Fusion}
Following Date et al.\cite{DJKMO}, we consider the fusion of the Vertex-Face correspondence relation \eqref{vertexface}. 

The fusion of the SOS weights is briefly summarized as follows. The SOS weight \eqref{face} satisfies 
\bea 
&&W\BW{a}{b}{d}{c}{0}=\delta_{b,d},\\
&&W\BW{a}{b}{d}{c}{-1}=0\qquad {\rm if}\ |a-c|=2,\\
&&W\BW{a}{a\pm 1}{a\pm 1}{a}{-1}=-W\BW{a}{a\pm 1}{a\mp 1}{a}{-1}.
\ena
Then if one defines
\bea
&&W_{21}\BW{a}{b}{d}{c}{u}=\sum_{d'}W\BW{a}{a'}{d}{d'}{u+1}W\BW{a'}{b}{d'}{c}{u}, \lb{fhfusion}
\ena
one can verify the following statements.

(i) The RHS of \eqref{fhfusion} is independent of the choice of $a'$ provided $|a-a'|=|a'-b|=1$.

(ii) For all $a, b, c, d$, 
\be
&&W_{21}\BW{a}{b}{d}{c}{-1}=0.
\en 
The $2\times 2$ fusion of the SOS weight is then given by the formula
\bea
&&W_{22}\BW{a}{b}{d}{c}{u}=\sum_{a'}W_{21}\BW{a}{b}{a'}{b'}{u-1}W_{21}\BW{a'}{b'}{d}{c}{u} 
\lb{fvhfusion}.
\ena
Here the RHS is independent of the choice of $b'$ provided $|b-b'|=|b'-c|=1$. Now the dynamical variables $a, b, c, d$ satisfies the extended admissible condition; $a_j-a_k\in \{2, 0, -2\}$ for any two adjacent local heights $a_j, a_k$. Furthermore the resultant SOS weight $W_{22}$ satisfies the face type YBE and defines the $2\times 2$ fusion SOS model. Explicit expressions of $W_{22}\BW{a}{b}{d}{c}{u}$ are given, for example, in \cite{KKW}. It satisfies the unitarity and crossing symmetry relations
\bea
&&\sum_{s}W_{22}\BW{a}{s}{d}{c}{-u} W_{22}\BW{a}{b}{s}{c}{u}=\delta_{b,d},\label{SOSunitarity}\\
&&W_{22}\BW{d}{c}{a}{b}{u}=\frac{(b,c)_2 \ g_a g_c}{(a,d)_2\ g_b g_d} \,
W_{22}\BW{a}{d}{b}{c}{-1-u}.
\label{SOScrossing}
\ena
Here $g_a=\vep_a\sqrt{[a]}$\, $\vep_a=\pm1,\ \vep_a\vep_{a+1}=(-)^a$ and 
\be
&&(a,b)_M=(b,a)_M=\left[\matrix{M\cr
                                 \frac{a-b+M}{2}\cr}\right]^{-1}
                                 \frac{\left[\frac{a+b-M}{2}, \frac{a+b+M}{2}
                                  \right]}{\sqrt{[a][b]}},\\
&&\left[\matrix{A\cr B\cr}\right]=\frac{[A][A-1]\cdots [A-B+1]}{[B][B-1]\cdots
[1]},\\
&&[A,B]=[A][A+1]\cdots [B] \quad (A<B),\qquad [A,A-1]=1.
\en 

The fusion of the intertwining vectors is given by
\bea
&&\psi^{(2)}(u)^a_b=\Pi\ \psi(u+1)^a_c\otimes \psi(u)^c_b \ \in V\otimes V.\lb{fusionpsi}
\ena
The RHS is independent of the choice of $c$ provided $|a-c|=|c-b|=1$.
Then using \eqref{vertexface}, \eqref{hfusion}, \eqref{vhfusion}, \eqref{fhfusion}, \eqref{fvhfusion} 
and \eqref{fusionpsi}, pne can show
\bea
&&\sum_{\mu_1',\mu_2'}
R^{(2,2)}(u-v)_{\mu_1 \mu_2}^{\mu_1' \mu_2'}\ 
\psi^{(2)}_{\mu_1'}(u)_{b}^{a}
\psi^{(2)}_{\mu_2'}(v)_{c}^{b}
=\sum_{b' \in {\mathbb{Z}}}
\psi^{(2)}_{\mu_2}(v)_{b'}^{a}
\psi^{(2)}_{\mu_1}(u)_{c}^{b'}
W_{22}\left(\left.
\begin{array}{cc}
a&b\\
b'&c
\end{array}\right|u-v\right).\nn\\
&&\lb{fusionvertexface}
\ena
Explicitly, the vector $\psi^{(2)}(u)^a_b$ is calculated as follows\cite{KKW}.
\begin{prop}
\begin{eqnarray*}
\left(\begin{array}{c}
\psi^{(2)}_2(u)_{n+2}^n\\
\psi^{(2)}_0(u)_{n+2}^n\\
\psi^{(2)}_{-2}(u)_{n+2}^n
\end{array}\right)=
\left(\begin{array}{c}
\vartheta_0\left(\left.
\frac{u-n+1
}{2r}
\right|\frac{\tau}{2}
\right)
\vartheta_0
\left(
\left.\frac{u-n-1}{2r}\right|\frac{\tau}{2}
\right)\\
{2}
\vartheta_0\left(\left.
\frac{u-n}{r}\right|
\tau\right)
\vartheta_0\left(
\left.\frac{1}{r}\right|\tau
\right)
\\
\vartheta_3\left(\left.
\frac{u-n+1
}{2r}
\right|\frac{\tau}{2}
\right)
\vartheta_3
\left(
\left.\frac{u-n-1}{2r}\right|\frac{\tau}{2}
\right)
\end{array}\right),
\end{eqnarray*}
\begin{eqnarray*}
\left(\begin{array}{c}
\psi^{(2)}_2(u)_{n}^n\\
\psi^{(2)}_0(u)_{n}^n\\
\psi^{(2)}_{-2}(u)_{n}^n
\end{array}\right)=
\left(\begin{array}{c}
\vartheta_0\left(\left.
\frac{u-n+1}{2r}\right|\frac{\tau}{2}
\right)
\vartheta_0\left(
\left.\frac{u+n+1}{2r}\right|\frac{\tau}{2}
\right)\\
{2}
\vartheta_0\left(\left.
\frac{n}{r}\right|
\tau
\right)
\vartheta_0
\left(
\left.\frac{u+1}{r}
\right|
\tau
\right)\\
\vartheta_3\left(\left.
\frac{u-n+1}{2r}\right|\frac{\tau}{2}
\right)
\vartheta_3\left(
\left.\frac{u+n+1}{2r}\right|\frac{\tau}{2}
\right)
\end{array}
\right),
\end{eqnarray*}
\begin{eqnarray*}
\left(\begin{array}{c}
\psi^{(2)}_2(u)_{n-2}^n\\
\psi^{(2)}_0(u)_{n-2}^n\\
\psi^{(2)}_{-2}(u)_{n-2}^n
\end{array}\right)=
\left(\begin{array}{c}
\vartheta_0\left(\left.
\frac{u+n+1}{2r}\right|\frac{\tau}{2}
\right)
\vartheta_0\left(
\left.\frac{u+n-1}{2r}\right|\frac{\tau}{2}
\right)
\\
{2}\vartheta_0\left(\left.
\frac{u+n}{r}\right|
\tau
\right)
\vartheta_0\left(
\left.\frac{1}{r}\right|
\tau
\right)
\\
\vartheta_3\left(\left.
\frac{u+n+1}{2r}\right|\frac{\tau}{2}
\right)
\vartheta_3\left(
\left.\frac{u+n-1}{2r}\right|\frac{\tau}{2}
\right)\end{array}\right).
\end{eqnarray*}
\end{prop}

\subsection{The dual intertwining vectors and their fusion}

Let us consider the dual vector $\psi^*(u)^a_b$ defined by 
\bea
&&\psi^*(u)^a_b\ v_{\vep}= \psi^*_\vep(u)^a_b,\qquad 
\psi^*_\vep(u)^a_b=-\vep\frac{a-b}{2[b][u]}C^2\ \psi_{-\vep}(u-1)^a_b \lb{dualintvec}
\ena
with $|a-b|=1$.
By a direct calculation, we verify the inversion relations
\bea
&&\sum_{\vep=\pm}\psi_\vep^*(u)^a_b\psi_\vep(u)^b_c=\delta_{a,c},\lb{inversiona}\\
&&\sum_{a=b\pm1}\psi_{\vep'}^*(u)^a_b\psi_\vep(u)^b_a=\delta_{\vep',\vep}.\lb{inversionb}
\ena
Hence we call $\psi^*(u)^a_b$  the dual intertwining vector. 
From the crossing symmetry properties of $R$ and $W$ the following vertex-face correspondence is held. 
\bea
&&\sum_{\vep_1',\vep_2'}
R(u-v)^{\vep_1 \vep_2}_{\vep_1' \vep_2'}\ 
\psi^*_{\vep_1'}(u)_{b}^{a}
\psi^*_{\vep_2'}(v)_{c}^{b}
=\sum_{s \in {\mathbb{Z}}}
\psi^*_{\vep_2}(v)_{b'}^{a}
\psi^*_{\vep_1}(u)_{c}^{b'}
W\left(\left.
\begin{array}{cc}
c&b'\\
b&a
\end{array}\right|u-v\right).\lb{vertexfacedual} 
\ena

The fusion of the dual intertwining vectors is given by\cite{KKW}
\begin{eqnarray}
\psi^{*(2)}(u)_a^b=
\sum_{c=a\pm 1} \psi^*(u+1)_a^c \otimes \psi^*(u)_c^b. \lb{fusiondualint} 
\end{eqnarray}
Then $\psi^{*(2)}(u)_a^b$ satisfies 
\bea
&&\Pi\ \psi^{*(2)}(u)_a^b=\psi^{*(2)}(u)_a^b\ \Pi. \lb{prdualint}
\ena
In the components, \eqref{fusiondualint} yields    
\begin{eqnarray}
\psi^{*(2)}_{\mu}(u)_a^b=
\sum_{c=a\pm 1} \psi^*_{\vep_1}(u+1)_a^c \psi^*_{\vep_2}(u)_c^b. \lb{compfusiondualint}
\end{eqnarray}
The relation \eqref{prdualint} indicates that the RHS of \eqref{compfusiondualint}is independent 
of the choice of $\vep_1, \vep_2$  proivided $\mu=\vep_1+\vep_2$. 
Then using \eqref{inversiona} and \eqref{inversionb}, it is easy to verify the following inversion relations.
\begin{prop}
\bea
&&\sum_{\vep=0,\pm2}\psi^{*(2)}_\vep(u)^a_b\psi^{(2)}_\vep(u)^b_c=\delta_{a,c},
\\
&&\sum_{a=b,b\pm 2}\psi^{*(2)}_{\vep'}(u)^a_b\psi^{(2)}_\vep(u)^b_a=\delta_{\vep', \vep}.
\ena
\end{prop}
Furthermore, in the similar way to the derivation of \eqref{fusionvertexface}, we obtain the fused form of \eqref{vertexfacedual}
\bea
&&\sum_{\mu_1',\mu_2'}
R^{(2,2)}(u-v)_{\mu_1 \mu_2}^{\mu_1' \mu_2'}\ 
\psi^{*(2)}_{\mu_1'}(u)_{b}^{a}
\psi^{*(2)}_{\mu_2'}(v)_{c}^{b}
=\sum_{b' \in {\mathbb{Z}}}
\psi^{*(2)}_{\mu_2}(v)_{b'}^{a}
\psi^{*(2)}_{\mu_1}(u)_{c}^{b'}
W_{22}\left(\left.
\begin{array}{cc}
c&b'\\
b&a
\end{array}\right|u-v\right).\nn\\
&&\lb{fusionvertexfacedual}
\ena 
The expression of $\psi^{*(2)}_\mu(u)_a^b\ (\mu=2, 0, -2)$ is evaluated as follows. 
\begin{prop}
\begin{eqnarray}
\left(\begin{array}{c}
\psi^{*(2)}_2(u)_{n+2}^n\\
\psi^{*(2)}_0(u)_{n+2}^n\\
\psi^{*(2)}_{-2}(u)_{n+2}^n
\end{array}\right)=\frac{C^4}{4[n+1][n+2][u][u+1]}
\left(\begin{array}{c}
\vartheta_3\left(\left.
\frac{u-n-1
}{2r}
\right|\frac{\tau}{2}
\right)^2\\
-\vartheta_3\left(\left.
\frac{u-n-1}{2r}\right|
\frac{\tau}{2}\right)
\vartheta_0\left(
\left.\frac{u-n-1}{2r}\right|\frac{\tau}{2}
\right)
\\
\vartheta_0\left(\left.
\frac{u-n-1
}{2r}
\right|\frac{\tau}{2}
\right)^2
\end{array}\right),\nn
\end{eqnarray}
\begin{eqnarray}
&&\left(\begin{array}{c}
\psi^{*(2)}_2(u)_{n}^n\\
\psi^{*(2)}_0(u)_{n}^n\\
\psi^{*(2)}_{-2}(u)_{n}^n
\end{array}\right)=-\frac{C^5}{4[n][n+1][n-1][u][u+1]}\nn\\
&&\qquad\times \left(\begin{array}{c}
\vartheta_3\left(\left.
\frac{u+n+1}{2r}\right|\frac{\tau}{2}
\right)
\vartheta_3\left(
\left.\frac{u-n-1}{2r}\right|\frac{\tau}{2}
\right)\vt{1}{n-1}{\tau}+\vartheta_3\left(\left.
\frac{u-n+1}{2r}\right|\frac{\tau}{2}
\right)
\vartheta_3\left(
\left.\frac{u+n-1}{2r}\right|\frac{\tau}{2}
\right)\vt{1}{n+1}{\tau}\\
-\vartheta_1\left(\left.
\frac{n}{r}\right|
\frac{\tau}{2}
\right)
\vartheta_2
\left(
\left.\frac{1}{r}
\right|
\frac{\tau}{2}
\right)
\vartheta_0
\left(
\left.\frac{u}{r}
\right|{\tau}
\right)\\
\vartheta_0\left(\left.
\frac{u+n+1}{2r}\right|\frac{\tau}{2}
\right)
\vartheta_0\left(
\left.\frac{u-n-1}{2r}\right|\frac{\tau}{2}
\right)\vt{1}{n-1}{\tau}+\vartheta_0\left(\left.
\frac{u-n+1}{2r}\right|\frac{\tau}{2}
\right)
\vartheta_0\left(
\left.\frac{u+n-1}{2r}\right|\frac{\tau}{2}
\right)\vt{1}{n+1}{\tau}
\end{array}
\right),\nn
\end{eqnarray}
\begin{eqnarray}
\left(\begin{array}{c}
\psi^{*(2)}_2(u)_{n-2}^n\\
\psi^{*(2)}_0(u)_{n-2}^n\\
\psi^{*(2)}_{-2}(u)_{n-2}^n
\end{array}\right)=\frac{C^4}{4[n-1][n-2][u][u+1]}
\left(\begin{array}{c}
\vartheta_3\left(\left.
\frac{u+n-1}{2r}\right|\frac{\tau}{2}
\right)^2
\\
-\vartheta_3\left(\left.
\frac{u+n-1}{2r}\right|
\frac{\tau}{2}
\right)
\vartheta_0\left(
\left.\frac{u+n-1}{2r}\right|
\frac{\tau}{2}
\right)
\\
\vartheta_0\left(\left.
\frac{u+n-1}{2r}\right|\frac{\tau}{2}
\right)^2
\end{array}\right).\nn
\end{eqnarray}
\end{prop}

Applying the crossing symmetry relations \eqref{crossing} and \eqref{SOScrossing} twice to \eqref{fusionvertexface}, we obtain the relation which should be compared with \eqref{fusionvertexfacedual}. Then fixing the suitable normalization function, we obtain\begin{thm}
\be
&&\psi^{*(2)}_{\vep}(u)^a_b=-\frac{C^4}{4[u][u+1]}\frac{\vtf{3}{0}{\tau}}{\vt{3}{1}{\tau}}\frac{g_a}{g_b (a,b)_2} \sum_{\vep'=0,\pm2}Q^{\vep'}_{\vep}\psi^{(2)}_{\vep'}(u-1)^a_b.
\en
\end{thm}

\vspace{3mm}
~\\
{\Large\bf Acknowledgements}~~

\noindent
The author would like to thank Michio Jimbo for sending his note and for stimulating discussions.
He is also grateful to Takeo Kojima and Robert Weston for their collaboration in the work \cite{KKW}.

\end{document}